\documentclass[11pt,a4paper]{article}
\usepackage[T1]{fontenc}
\usepackage{lmodern}
\usepackage{amsmath,amssymb,amsthm}
\usepackage[margin=2.9cm]{geometry}
\usepackage[hidelinks]{hyperref}
\hypersetup{pdftitle={Co-amenability and local reflexivity of C*-algebras},
  pdfauthor={Yafei Zhao and Zhe Dong}}
\numberwithin{equation}{section}

\newtheorem{theorem}{Theorem}[section]
\newtheorem{proposition}[theorem]{Proposition}
\newtheorem{lemma}[theorem]{Lemma}
\newtheorem{corollary}[theorem]{Corollary}
\theoremstyle{definition}
\newtheorem{definition}[theorem]{Definition}

\newcommand{\C}{\mathbb C}
\newcommand{\Z}{\mathbb Z}
\newcommand{\BH}{\mathcal B}
\newcommand{\OTwo}{\mathcal O_2}
\newcommand{\op}{\mathrm{op}}
\newcommand{\id}{\mathrm{id}}
\newcommand{\tmin}{\otimes_{\min}}
\newcommand{\ucp}{\mathrm{u.c.p.}}

\title{Co-amenability and local reflexivity of $C^*$-algebras}
\author{Yafei Zhao and Zhe Dong\thanks{Corresponding author.}}
\date{}

\begin{document}
\maketitle

\begin{abstract}
Let $A\subset B$ be a unital inclusion of $C^*$-algebras. We prove that
local reflexivity passes from $A$ to $B$ if the inclusion
$A^{**}\subset B^{**}$ is co-amenable. It follows that local
reflexivity is preserved by minimal tensor products with exact
$C^*$-algebras, by inductive limits over directed sets, and by
reduced crossed products by discrete amenable groups.
We also prove that reduced free products with respect to states with faithful
GNS representations preserve local reflexivity and that every locally
reflexive $C^*$-algebra embeds into a simple locally reflexive $C^*$-algebra.
These results answer Questions~(b), (c), (d), (f) and~(h) in Kirchberg's
2011 list.
\end{abstract}



\section{Introduction}

Local reflexivity is a finite-dimensional approximation property of
operator spaces. For a unital $C^*$-algebra $A$, it can be formulated
as the existence of unital completely positive approximations of the
inclusions of finite-dimensional operator systems in $A^{**}$ by
maps taking values in $A$. The convergence is pointwise in
$\sigma(A^{**},A^*)$ \cite[Theorem~5.1]{EH}. In his 2011 Copenhagen
lecture, Kirchberg asked whether local reflexivity is preserved by
minimal tensoring with the CAR algebra $M_{2^\infty}$, by taking
inductive limits, and by crossing with $\Z$ under an arbitrary automorphism
\cite[slide~12, Questions~(b), (c), and~(d)]{Kir2011}. He also stated the
equivalent formulation of Question~(b) with an arbitrary exact tensor
factor.

Our main result is a permanence theorem for inclusions. An inclusion
of von Neumann algebras $R\subset M$ is \emph{co-amenable} if in the
standard representation of $M$ on $L^2(M)$ there is a conditional
expectation from $R'$ onto $M'$. We allow nonnormal conditional
expectations.

\begin{theorem}[Co-amenability transfer]\label{thm:transfer}
Let $A\subset B$ be a unital inclusion of $C^*$-algebras. Identify
$A^{**}$ with its image under the normal extension of the inclusion in
$B^{**}$. If $A^{**}\subset B^{**}$ is co-amenable and $A$ is locally
reflexive, then $B$ is locally reflexive.
\end{theorem}

\begin{corollary}[Kirchberg's Questions (b), (c), and (d)]\label{cor:kirchberg}
Let $A$ be a $C^*$-algebra.
\begin{enumerate}
\item[\rm(b)] If $A$ is locally reflexive and $D$ is an exact $C^*$-algebra, then
  $A\tmin D$ is locally reflexive. In particular, this holds for
  $D=M_{2^\infty}$.
\item[\rm(c)] If $A$ is the inductive limit of a directed system of locally
  reflexive $C^*$-algebras, then $A$ is locally reflexive.
\item[\rm(d)] If $A$ is locally reflexive and a discrete amenable group $G$ acts on
  $A$ by automorphisms, then $A\rtimes_{\mathrm r}G$ is locally reflexive.
  In particular, $A\rtimes_\alpha\Z$ is locally reflexive for every
  automorphism $\alpha$ of $A$.
\end{enumerate}
\end{corollary}

For~(c) we prove more generally that a $C^*$-algebra is
locally reflexive if each finite subset can be approximated arbitrarily
well in norm within a locally reflexive $C^*$-subalgebra.
The subalgebras used for different finite sets need not be nested.

Kirchberg also asked about reduced free products and embeddings into
simple algebras. His lecture lists the implications
$\mathrm{(d)}\Rightarrow\mathrm{(f)}\Rightarrow\mathrm{(h)}$
\cite[slides~12--13]{Kir2011}. In Section~\ref{sec:free-simple}
we derive both conclusions from Corollary~\ref{cor:kirchberg}
using the Toeplitz--Pimsner construction of Dykema and Shlyakhtenko
\cite{DS}.

\begin{corollary}[Kirchberg's Questions (f) and (h)]\label{cor:kirchberg-fh}
\leavevmode
\begin{enumerate}
\item[\rm(f)] Let $(A_i,\varphi_i)_{i\in I}$ be a family of unital locally
  reflexive $C^*$-algebras with states whose GNS representations are
  faithful. Their reduced free product over $\C$ is locally reflexive.
\item[\rm(h)] Every locally reflexive $C^*$-algebra embeds into a simple
  locally reflexive $C^*$-algebra.
\end{enumerate}
\end{corollary}

Dong and Dong constructed non-locally
reflexive reduced group $C^*$-algebras and proved a local reflexivity
characterization for a family of wreath products
\cite[Theorem~5.3 and Corollary~5.4]{DD}, thereby answering
Question~(a) of Kirchberg's list in part.

The proof of Theorem~\ref{thm:transfer} uses the weak containment
criterion of Bannon, Marrakchi, and Ozawa \cite[Corollary~A.3]{BMO}
and the coefficient approximation theorem of Anantharaman--Delaroche
\cite[Theorem~2.2 and \S2.7]{AD}. For a co-amenable inclusion
$A^{**}\subset B^{**}$, these results approximate the identity on
$B^{**}$ point ultraweakly by completely positive maps that factor
through matrix algebras over $A^{**}$.

Hasegawa obtained conditional expectations between commutants in his
work on relative nuclearity, with applications to nuclear tensor
factors and crossed products by discrete amenable groups \cite{Has}.
His relative completely positive approximations converge pointwise in
norm on $B$. Local reflexivity requires approximating the inclusion
of each finite-dimensional operator system $F\subset B^{**}$ by maps
into $B$, with convergence tested against $B^*$. Starting with the
coefficient maps on $B^{**}$, we use Kaplansky density and local
reflexivity of $A$ to obtain unital completely positive maps $F\to B$
with this convergence.

\section{Local reflexivity of \texorpdfstring{$C^*$}{C*}-algebras}

All unital inclusions in this paper preserve the unit. A conditional
expectation onto a unital subalgebra is a unital completely positive
projection onto that subalgebra. Its bimodule property follows from
Tomiyama's theorem \cite{Tom}. We use $M_n(A)$ and $M_{n,1}(A)$ for
matrices and columns over $A$, respectively.

\begin{definition}
A $C^*$-algebra $A$ is \emph{locally reflexive} if every complete
contraction from a finite-dimensional operator space into $A^{**}$
is the pointwise $\sigma(A^{**},A^*)$ limit of complete contractions
with values in $A$.
\end{definition}

By Effros and Haagerup \cite[Theorem~5.1]{EH}, a unital $A$ is locally
reflexive precisely when every finite-dimensional operator system
$F\subset A^{**}$ admits $\ucp$ maps into $A$ converging point
ultraweakly to the inclusion. For arbitrary $A$ (see
\cite[\S5, pp.~119--120]{EH} for the unitization argument), their property
$(C'')$ characterizes local reflexivity by contractivity of the canonical map
\begin{equation}\label{eq:property-c}
  j_{A,D}\colon A^{**}\odot D \longrightarrow (A\tmin D)^{**}
\end{equation}
for the minimal norm and every $C^*$-algebra $D$. Here $j_{A,D}$
extends $a\otimes d$ normally in the first variable. A $C^*$-algebra is
locally reflexive if all its separable $C^*$-subalgebras are locally
reflexive \cite[slide~11]{Kir2011}.

\begin{lemma}\label{lem:permanence}
Local reflexivity passes to $C^*$-subalgebras, quotients, and unitizations.
A unital extension with locally reflexive ideal and quotient is locally
reflexive if its quotient map has a unital $*$-homomorphic section.
If $A$ is unital and locally reflexive, then for every
finite-dimensional operator system $F\subset M_n(A^{**})$ its
inclusion is the point-ultraweak limit of $\ucp$ maps
$F\to M_n(A)$.
\end{lemma}

\begin{proof}
For $A_0\subset A$, both $A_0\hookrightarrow A$ and
$A_0\tmin D\hookrightarrow A\tmin D$ induce isometric bidual
inclusions. By naturality, $j_{A_0,D}$ is the restriction of $j_{A,D}$
under these inclusions. Injectivity of the minimal norm then implies
property $(C'')$ for $A_0$ whenever it holds for $A$.

For a quotient map $q\colon A\to Q$, the normal map $q^{**}$ identifies
a central summand of $A^{**}$ with $Q^{**}$. Lift a complete contraction
into $Q^{**}$ to this summand, apply local reflexivity of $A$, and
compose the approximating maps with $q$.

A $*$-split extension $0\to I\to E\to Q\to 0$ remains split exact
after minimal tensoring with any $C^*$-algebra. If $I$ and $Q$ are
locally reflexive, then $E$ has property $(C'')$ by the criterion of
Effros and Haagerup \cite[Proposition~5.5]{EH}. The unitization
assertion follows by applying this to $0\to I\to\widetilde I\to\C\to0$.

Finally, let $F_0\subset A^{**}$ be the finite-dimensional operator
system spanned by $1_A$ and all entries of matrices in $F$. Take
$\ucp$ maps $\theta_i\colon F_0\to A$ converging point ultraweakly to
the inclusion. The restrictions of $\id_{M_n}\otimes\theta_i$ to
$F$ are $\ucp$ and converge entrywise ultraweakly. Every normal
functional on $M_n(A^{**})$ is a finite sum of normal functionals
on the entries, so the convergence also holds on $F$ in the
point-ultraweak topology.
\end{proof}

\section{Proof of the transfer theorem}

\subsection{Completely positive coefficients}

\begin{lemma}\label{lem:coefficients}
Let $R\subset M$ be a unital co-amenable inclusion of von Neumann
algebras. There is a net of maps $\Phi_i\colon M\to M$ converging
point ultraweakly to $\id_M$, each of the form
\begin{equation}\label{eq:normalized-coefficient}
  \begin{gathered}
    \Phi_i(x)=v_i^*\alpha_i(x)v_i,\\
    \alpha_i\colon M\to M_{n_i}(R)\ \text{normal and } \ucp,
    \qquad v_i\in M_{n_i,1}(M),\quad \|v_i\|\leq 1 .
  \end{gathered}
\end{equation}
\end{lemma}

\begin{proof}
Let $H={}_M L^2(M)_R$ be the correspondence with the standard
left action of $M$ and the restricted right action of $R$. If
$J_M$ denotes the standard conjugation \cite{Haa}, the right
action is $\xi r=J_M r^*J_M\xi$. Its commutant is
\[
  Q=\BH_{R^{\op}}(H)=J_MR'J_M .
\]
By co-amenability, we can choose an expectation $P\colon R'\to M'$.
Then $y\mapsto J_M P(J_MyJ_M) J_M$ defines a conditional expectation
$Q\to M$. Conjugation by $J_M$ is conjugate-linear and preserves
products and positivity, also entrywise at each matrix level.
Hence this composite of $P$ with two conjugations is a linear
completely positive projection. By \cite[Corollary~A.3]{BMO},
\begin{equation}\label{eq:weak-containment}
  {}_M L^2(M)_M \prec {}_M L^2(Q)_M .
\end{equation}
Here $\prec$ denotes weak containment of correspondences. The
corollary permits a nonnormal expectation. By the standard-form
identification recalled in \cite[\S2.1]{AD}, $L^2(Q)$ as an $M$-$M$
correspondence is $H\mathbin{\overline{\otimes}_R}\overline H$.

Apply \cite[Theorem~2.2 and \S2.7]{AD} to
$H={}_M L^2(M)_R$. In the notation of \S2.7 of that paper,
$S_j$ commutes with the left $M$-action on $L^2(M)$, so
$w_j=J_MS_jJ_M\in M$. The other intertwiners
$T_j\colon L^2(R)\to H$ are right $R$-linear, so $T_j^*xT_k\in R$.
With these identifications the coefficient in \S2.7 is
\begin{equation}\label{eq:raw-coefficient}
  x\longmapsto\sum_{j,k=1}^n w_j^*(T_j^*xT_k)w_k .
\end{equation}
By \eqref{eq:weak-containment}, the identity on $M$ therefore lies in
the point-ultraweak closure of the cone of finite sums of
\eqref{eq:raw-coefficient}. The cone is closed under
$\Phi\mapsto b^*\Phi(\,\cdot\,)b$ for $b\in M$, since this replaces
each $w_j$ by $w_jb$. It follows from \cite[Lemma~1.6]{AD}
that the approximants $\Phi$ can be chosen with $\Phi(1)\leq 1$.

Consider one coefficient in the finite sum for such a $\Phi$. With
the row operator $T$ and the column $w$ given by
\[
  T=[T_1\ \cdots\ T_n]\colon L^2(R)^{\oplus n}\to H,
  \qquad w=(w_j)_{j=1}^n,
\]
the coefficient at $x$ is $w^*(T^*xT)w$.
Let $T=Vh^{1/2}$ be the polar decomposition, where
$h=T^*T\in M_n(R)$. The partial isometry $V$ remains a right
$R$-module intertwiner. Indeed, for each unitary of the right
action the intertwining equation passes to $V$ by uniqueness
of the polar decomposition, as $h$ commutes with that action.
Thus $V^*xV\in M_n(R)$ and $p=V^*V$ is the support projection of
$h$ in $M_n(R)$. We choose a normal state $\omega$ on $M$. The map
\[
  \alpha(x)=V^*xV+\omega(x)(1-p)
\]
is normal, unital, and completely positive. If $v=h^{1/2}w$,
with products computed using $R\subset M$, then $(1-p)v=0$ and
\[
  v^*\alpha(x)v=w^*(T^*xT)w,\qquad
  v^*v=w^*(T^*T)w .
\]
For the finite sum, use the block diagonal of the individual
maps $\alpha$ and the column obtained by stacking the individual
columns $v$. This is the factorization in
\eqref{eq:normalized-coefficient}, with
$v^*v=\Phi(1)\leq 1$, so the stacked column is a contraction.
Applying the construction to each approximant proves the lemma.
\end{proof}

\subsection{Approximations with values in the original algebra}

\begin{theorem}[Theorem~\ref{thm:transfer}]\label{thm:transfer-main}
Let $A\subset B$ be a unital inclusion of $C^*$-algebras. Identify
$A^{**}$ with its image under the normal extension of the inclusion in
$B^{**}$. If $A^{**}\subset B^{**}$ is co-amenable and $A$ is locally
reflexive, then $B$ is locally reflexive.
\end{theorem}

\begin{proof}
Throughout the proof, $M=B^{**}$ and $R=A^{**}$, the latter embedded
normally in $M$. Let $F\subset M$ be a finite-dimensional
operator system. We first approximate its inclusion by completely positive
contractions $F\to B$. Consider finite test sets $\mathcal X\subset F$ and
$\mathcal L\subset M_*=B^*$ and a tolerance $\varepsilon>0$. By
Lemma~\ref{lem:coefficients}, we can choose $x\mapsto v^*\alpha(x)v$ within
$\varepsilon/3$ of the inclusion on these tests.

By Kaplansky density \cite{Kap}, some net of contractions
$c\in M_{n,1}(B)$ converges strong$^*$ to $v$: take the upper-right
corners of self-adjoint contractions in $M_{n+1}(B)$ approximating
$\left(\begin{smallmatrix}0&v\\v^*&0\end{smallmatrix}\right)$.
Multiplication is jointly strong$^*$ continuous on bounded sets.
We may thus choose $c$ so that
$|\ell(c^*\alpha(x)c-v^*\alpha(x)v)|<\varepsilon/3$ on the tests.

By Lemma~\ref{lem:permanence}, the inclusion of the finite-dimensional
operator system $\alpha(F)$ in $M_n(R)$ is the point-ultraweak limit of
$\ucp$ maps into $M_n(A)$. Since $y\mapsto\ell(c^*yc)$ is normal on
$M_n(R)$, we can compose one such map with $\alpha|_F$ to obtain a $\ucp$ map
$\psi\colon F\to M_n(A)$ with
$|\ell(c^*(\psi(x)-\alpha(x))c)|<\varepsilon/3$ on the tests.

Then $x\mapsto c^*\psi(x)c$ is a completely positive contraction into
$B$ and differs from the inclusion by less than $\varepsilon$ on the
tests. Letting the tests and the tolerance vary, we obtain a net
$(\gamma_i)$ of these contractions converging point ultraweakly to the
inclusion. In particular, $\gamma_i(1)\to1$ ultraweakly. For any state $\omega$ on $F$,
\begin{equation}\label{eq:output-cp}
  \theta_i(x)=\gamma_i(x)+\omega(x)(1-\gamma_i(1)),\qquad x\in F,
\end{equation}
defines a $\ucp$ map $F\to B$ with the same point-ultraweak limit.
\end{proof}

\section{Tensor products, inductive limits, and crossed products}

Anantharaman--Delaroche studied relative amenability for von Neumann
crossed products and spatial tensor products \cite[Proposition~3.4 and
\S4.2]{AD}. In Propositions~\ref{prop:nuclear} and~\ref{prop:crossed}
we consider the canonical inclusions of $A^{**}$ into the full
biduals of the corresponding $C^*$-algebraic tensor product and
reduced crossed product.

\subsection{Exact tensor factors}

\begin{proposition}\label{prop:nuclear}
If $A$ and $D$ are unital $C^*$-algebras and $D$ is nuclear,
the inclusion $A^{**}\subset (A\tmin D)^{**}$ induced by
$a\mapsto a\otimes 1_D$ is co-amenable.
\end{proposition}

\begin{proof}
Let $B=A\tmin D$ and $M=B^{**}$. In the standard representation
of $M$, denote by $R$ and $N$ the images of $A^{**}$ and $D^{**}$
under the normal extensions of the embeddings of the two
factors. The images commute by separate ultraweak continuity of
multiplication, and they generate $M$ as a von Neumann algebra
since the two original factors generate $B$. In particular,
$M'=R'\cap N'$.

The bidual of the embedding $d\mapsto 1_A\otimes d$ is faithful, so
$N$ is normally $*$-isomorphic to $D^{**}$. Thus $N$ is injective by
nuclearity of $D$ \cite{CE}, and its commutant $N'$ is injective as
well \cite{EL}. Hence there is a conditional expectation
$P\colon\BH(L^2(M))\to N'$. For $x\in R'$ and $r\in R$ we have
$r\in N'$ and hence, by the $N'$-bimodule property of $P$,
\[
  rP(x)=P(rx)=P(xr)=P(x)r .
\]
Thus $P(R')\subset R'\cap N'=M'$. The restriction of $P$ to
$R'$ fixes $M'$ pointwise and is the required expectation.
\end{proof}

Hasegawa's tensor product example \cite[Example~5.2.1]{Has} uses a state
on $D$ with faithful GNS representation. In Proposition~\ref{prop:nuclear}
one may take $D=\widetilde{c_0(I)}$ for uncountable $I$, where the
GNS representation of each state annihilates all but countably many of
the central coordinate projections.

\begin{corollary}[Corollary~\ref{cor:kirchberg}(b)]\label{cor:tensor}
If $A$ is locally reflexive and $D$ is an exact $C^*$-algebra, then
$A\tmin D$ is locally reflexive. In particular, this holds for
$D=M_{2^\infty}$.
\end{corollary}

\begin{proof}
We first treat nuclear $D$. For unital $A$ and unital $D$, apply
Theorem~\ref{thm:transfer-main} and Proposition~\ref{prop:nuclear}.
In general, adjoin units as needed. Unitization preserves local
reflexivity by Lemma~\ref{lem:permanence} and preserves nuclearity.
The original tensor product is a $C^*$-subalgebra of the tensor
product of the unital algebras, so it is locally reflexive by
Lemma~\ref{lem:permanence}. The CAR algebra $M_{2^\infty}$ is
nuclear, so this case includes the original formulation of
Question~(b).

Now let $D$ be exact. If $D$ is separable, then $D$ embeds into the
nuclear algebra $\OTwo$ by the Kirchberg--Phillips embedding theorem
\cite{KP}, and the conclusion follows from the inclusion
$A\tmin D\subset A\tmin\OTwo$ and hereditary local reflexivity. For
arbitrary $D$, let $C$ be any separable $C^*$-subalgebra of
$A\tmin D$. Approximate the elements of a countable dense set of $C$
by sequences of finite sums of elementary tensors. All the
$D$-coefficients in these sums belong to a separable $C^*$-subalgebra
$D_0\subset D$. Exactness passes to $D_0$, and $C$ is contained in
the closed subalgebra $A\tmin D_0$. By the separable case,
$A\tmin D_0$ is locally reflexive, and hence so is $C$. Separable
determination now shows that $A\tmin D$ is locally reflexive.
\end{proof}

\subsection{Local approximation and inductive limits}\label{sec:inductive}

\begin{lemma}[Diagonal inclusion]\label{lem:diagonal}
Let $(A_i)_{i\in I}$ be $C^*$-subalgebras of $A$, indexed by an infinite
set $I$. Suppose an ultrafilter $\mathcal U$ on $I$ satisfies
$\lim_{i\to\mathcal U}\operatorname{dist}(a,A_i)=0$ for each $a\in A$.
Consider the $C^*$-algebras
\[
  D_0=\bigoplus_{i\in I}^{c_0}A_i,\qquad
  B_0=\mathcal K(\ell^2(I))\tmin A .
\]
Using the standard matrix units on $\ell^2(I)$, embed $D_0$ in $B_0$
by $(a_i)\mapsto\sum_i e_{ii}\otimes a_i$. Let $B$ be the unitization
of $B_0$ and $D=D_0+\C 1_B$. Then $D^{**}\subset B^{**}$ is co-amenable.
\end{lemma}

\begin{proof}
We may assume $A\ne 0$. The diagonal sum converges in norm and has
norm $\sup_i\|a_i\|$. Represent $B^{**}$ in standard form on $H$, and denote by $\rho$ the
restriction of this representation to $B$. Let $H_1$ be the essential
subspace $\overline{\rho(B_0)H}$ of $\rho|_{B_0}$, and let
$H_0=H_1^\perp$. On $H_0$ the representation factors through the
scalar quotient. On $H_1$ we may identify $\rho|_{B_0}$ with
$\id_{\mathcal K}\otimes\pi$ on $\ell^2(I)\otimes L$, where $\pi$ is a
nondegenerate representation of $A$ on $L$. Indeed, the multiplier
matrix units identify all coordinate subspaces with $L$, and the
multiplier copy of $A$ commutes with these matrix units. The support
projection of $\rho(B_0)$ belongs to $\rho(B)''$, so
\[
  \rho(B)'=(1_{\ell^2(I)}\otimes\pi(A)')\oplus\BH(H_0).
\]

For $i\in I$, let $V_i\colon L\to H_1\subset H$ insert a vector into
coordinate $i$. For $x\in\rho(D)'$, define $x_i=V_i^*xV_i$. Compressing
the commutation relation between $x$ and $\rho(e_{ii}\otimes b)$ to
coordinate $i$ shows that $[x_i,\pi(b)]=0$ for every $b\in A_i$. Thus
for each $a\in A$,
\[
  \|[x_i,\pi(a)]\|\leq 2\|x\|\operatorname{dist}(a,A_i).
\]
The ball of radius $\|x\|$ in $\BH(L)$ is ultraweakly compact. Let
$y(x)$ be the ultraweak limit of $x_i$ along $\mathcal U$. By the
estimate and separate ultraweak continuity of multiplication,
$y(x)\in\pi(A)'$. If $P_0$ is the projection onto $H_0$, define
\[
  E(x)=(1_{\ell^2(I)}\otimes y(x))\oplus P_0x|_{H_0}.
\]
Compression is completely positive, and positivity at every matrix
level passes to the ultraweak ultralimit. Thus $E$ is $\ucp$ and
takes values in $\rho(B)'$.

For $x=(1\otimes z)\oplus w\in\rho(B)'$ we have $x_i=z$ for all $i$
and $P_0x|_{H_0}=w$, so $E(x)=x$. Hence $E$ is a conditional
expectation from $\rho(D)'$ onto $\rho(B)'$. These are the commutants
of $D^{**}$ and $B^{**}$ in the standard representation of $B^{**}$,
as required.
\end{proof}

\begin{theorem}[Local approximation]\label{thm:local-approximation}
Suppose that for every finite set $S\subset A$ and every $\varepsilon>0$
there is a locally reflexive $C^*$-subalgebra $C\subset A$ with
$\operatorname{dist}(a,C)<\varepsilon$ for all $a\in S$.
Then $A$ is locally reflexive.
\end{theorem}

\begin{proof}
For each $i=(S,n)$, with $S\subset A$ finite and $n\geq 1$, we choose
such a subalgebra $A_i$ with tolerance $1/n$. Order the pairs by
inclusion of the finite sets and by the usual order on the integers.
An ultrafilter containing every tail of this directed set satisfies
the hypothesis of Lemma~\ref{lem:diagonal}.

We first check local reflexivity of $D_0$ directly. Its bidual is
$\prod_i A_i^{**}$ with predual $\bigoplus_i^{\ell^1} A_i^*$.
Let $u\colon F\to D_0^{**}$ be a complete contraction on a
finite-dimensional operator space. For a finite set of tests from
$D_0^*$, we choose a common finite set of coordinates $J$ so that each
test has arbitrarily small $\ell^1$-norm outside $J$.
For each $i\in J$, local reflexivity of $A_i$ provides complete
contractions $F\to A_i$ approximating the $i$th coordinate of $u$
pointwise ultraweakly. Combining these maps on $J$ and setting the
other coordinates to zero produces a complete contraction into $D_0$.
For $t\in F$ and $(f_i)\in D_0^*$, the discarded coordinates contribute
at most $\|t\|\sum_{i\notin J}\|f_i\|$. The finite-coordinate
approximations therefore converge pointwise in the weak-star topology,
as required.

It follows by unitization that $D$ is locally reflexive
(Lemma~\ref{lem:permanence}). By Lemma~\ref{lem:diagonal} and
Theorem~\ref{thm:transfer-main}, $B$ is locally reflexive.
For any $j\in I$, the map $a\mapsto e_{jj}\otimes a$ embeds $A$ in $B$.
Local reflexivity passes to this copy of $A$ by
Lemma~\ref{lem:permanence}.
\end{proof}

\begin{corollary}[Corollary~\ref{cor:kirchberg}(c)]\label{cor:inductive}
If $A$ is the inductive limit of a directed system of locally
reflexive $C^*$-algebras, then $A$ is locally reflexive.
\end{corollary}

\begin{proof}
Let $\iota_i\colon A_i\to A$ be the canonical maps into the inductive
limit. Each $\iota_i(A_i)$ is a quotient of $A_i$ and is locally
reflexive by Lemma~\ref{lem:permanence}. These images form a directed
family with dense union, so $A$ is locally reflexive by
Theorem~\ref{thm:local-approximation}.
\end{proof}

\subsection{Crossed products by discrete amenable groups}

\begin{proposition}\label{prop:crossed}
Let a discrete amenable group $G$ act by automorphisms on a
unital $C^*$-algebra $A$. The inclusion
$A^{**}\subset (A\rtimes_{\mathrm r} G)^{**}$ is co-amenable.
\end{proposition}

\begin{proof}
Let $B=A\rtimes_{\mathrm r} G$, $R=A^{**}$ and $M=B^{**}$.
In the standard representation of $M$, let $u_g\in B$ be the
canonical unitaries. By the covariance relations,
$u_g R u_g^*=R$. Moreover, $M$ is generated as a von Neumann
algebra by $R$ and the $u_g$. It follows that
\begin{equation}\label{eq:fixed}
  (R')^{\operatorname{Ad}u(G)}
  =R'\cap\{u_g:g\in G\}'=M' .
\end{equation}
Take a left F{\o}lner net $(F_i)_i$ of nonempty finite subsets
of $G$. Each map
\[
  P_i(x)=\frac{1}{|F_i|}\sum_{g\in F_i}u_g x u_g^*,
  \qquad x\in R',
\]
is $\ucp$ from $R'$ to $R'$ and fixes $M'$. For $h\in G$,
\[
  \|u_h P_i(x)u_h^*-P_i(x)\|
    \leq \frac{|hF_i\mathbin{\triangle}F_i|}{|F_i|}\|x\|
    \longrightarrow 0 .
\]
Compactness of the product of the ultraweak balls
$\{y\in R':\|y\|\leq\|x\|\}$ over $x\in R'$ gives a point-ultraweak
cluster point $P$ of the maps $P_i$. It is $\ucp$
and fixes $M'$ pointwise. The estimate gives $u_hP(x)u_h^*=P(x)$ for
every $h\in G$, so $P$ is a conditional expectation onto $M'$
by \eqref{eq:fixed}.
\end{proof}

\begin{corollary}[Corollary~\ref{cor:kirchberg}(d)]\label{cor:crossed}
If $A$ is locally reflexive and a discrete amenable group $G$ acts on
$A$ by automorphisms, then $A\rtimes_{\mathrm r}G$ is locally reflexive.
In particular, $A\rtimes_\alpha\Z$ is locally reflexive for every
automorphism $\alpha$ of $A$.
\end{corollary}

\begin{proof}
For unital $A$, combine Proposition~\ref{prop:crossed} with
Theorem~\ref{thm:transfer-main}. For nonunital $A$, extend the action to
$\widetilde A$ by fixing the unit. A faithful nondegenerate
representation of $A$ extends faithfully to $\widetilde A$ on the same
space. The associated regular representations agree on finite sums with
coefficients in $A$. Hence
$A\rtimes_{\rm r}G\subset\widetilde A\rtimes_{\rm r}G$ isometrically.
The assertion follows from the unital case by Lemma~\ref{lem:permanence}.

For $G=\Z$ the full and reduced crossed products coincide by
amenability, so the assertion for $A\rtimes_\alpha\Z$ follows.
\end{proof}

\section{Reduced free products and simple embeddings}\label{sec:free-simple}

For a Hilbert $D$-bimodule $X$, let $\mathcal T(X)$ denote its
Toeplitz--Pimsner algebra, generated by the copy of $D$ and the creation
operators on the full Fock module. This is the algebra denoted by $E(X)$
in \cite{DS}.

\begin{proposition}\label{prop:toeplitz}
Let $D$ be a unital locally reflexive $C^*$-algebra and let $X$ be a
full right Hilbert $D$-module with a faithful unital left action of $D$
by adjointable operators. Then $\mathcal T(X)$ is locally reflexive.
\end{proposition}

\begin{proof}
We follow the proof of Dykema and Shlyakhtenko's exactness theorem
\cite[Theorem~3.1]{DS}, using their filtration of the fixed-point algebra
and their crossed-product description. For local reflexivity the required
permanence properties follow from Lemma~\ref{lem:permanence} and
Corollaries~\ref{cor:tensor}, \ref{cor:inductive}, and~\ref{cor:crossed}.

We first show that $\mathcal K(Z)$ is locally reflexive for every right
Hilbert $D$-module $Z$. For each finite subset $S\subset Z$, let $Z_S$
be its closed right $D$-span. By Kasparov stabilization
\cite[Theorem~2]{Kas}, $\mathcal K(Z_S)$ embeds in
$\mathcal K(\ell^2)\tmin D$, and hence is locally reflexive by
Corollary~\ref{cor:tensor} and Lemma~\ref{lem:permanence}. Through the
rank-one formulas, the algebras $\mathcal K(Z_S)$ form a directed
family of subalgebras of $\mathcal K(Z)$ whose union contains all
finite-rank operators. Hence $\mathcal K(Z)$ is locally reflexive by
Corollary~\ref{cor:inductive}.

Next let $Y=X\oplus D$. The inclusion $X\subset Y$ induces an
embedding $\mathcal T(X)\subset\mathcal T(Y)$. Let $C$ be the
fixed-point algebra of the gauge action on $\mathcal T(Y)$, and for
$n\geq0$ let $C_n$ be the closed span of $D$ and the words with $k$
creators followed by $k$ adjoints, $1\leq k\leq n$. Then $C_0=D$, $C$
is the closure of the increasing union of the $C_n$, and for $n\geq1$
there is a split exact sequence \cite[proof of Claim~3.5]{DS}
\[
  0\longrightarrow \mathcal K(Y^{\otimes_D n})
    \longrightarrow C_n\longrightarrow C_{n-1}\longrightarrow 0 .
\]
The quotient map is compression to Fock levels $0,\ldots,n-1$. It is
faithful on $C_{n-1}$, since on each higher level the action of
$C_{n-1}$ is its action on level $n-1$ tensored with the identity.
Its kernel is the closed span of the words with $n$ creators, which
restriction to level $n$ identifies with $\mathcal K(Y^{\otimes_D n})$.
The inclusion of $C_{n-1}$ splits the sequence. Since $C_0=D$ and the
kernels are locally reflexive by the preceding paragraph,
Lemma~\ref{lem:permanence} and induction on $n$ show that each $C_n$
is locally reflexive. Hence
$C$ is locally reflexive by Corollary~\ref{cor:inductive}.

Finally, let $L$ be creation by $0\oplus1_D$, and define
$\beta(c)=LcL^*$ for $c\in C$. By \cite[Claim~3.4 and \S2]{DS},
$\beta$ is an injective endomorphism of $C$, and $\mathcal T(Y)$ is a
corner of $C_\infty\rtimes_{\bar\beta}\Z$, where
$C_\infty=\varinjlim(C,\beta)$ and $\bar\beta$ is the induced
automorphism. The algebra $C_\infty$ is locally reflexive by
Corollary~\ref{cor:inductive}, and so is its crossed product by
Corollary~\ref{cor:crossed}. Local reflexivity passes to the corner
$\mathcal T(Y)$ and then to $\mathcal T(X)$ by
Lemma~\ref{lem:permanence}.
\end{proof}

\begin{corollary}[Corollary~\ref{cor:kirchberg-fh}(f)]\label{cor:free}
Let $(A_i,\varphi_i)_{i\in I}$ be a family of unital locally
reflexive $C^*$-algebras with states whose GNS representations are
faithful. Their reduced free product over $\C$ is locally reflexive.
\end{corollary}

\begin{proof}
We begin with two factors. Let $D=A_1\oplus A_2$, which is locally
reflexive by coordinatewise approximation. By
\cite[Proposition~5.1]{DS}, applied with amalgamation over $\C$, there
is a Hilbert $D$-bimodule $X$ such that the reduced free product of $(A_1,\varphi_1)$ and
$(A_2,\varphi_2)$ is a quotient of a $C^*$-subalgebra of
$\mathcal T(X)$. In that construction the vector represented by
$1_D\otimes 1_D$ has inner product $1_D$ with itself, so $X$ is full.
Writing $(\pi_i,H_i)$ for the GNS representation of $\varphi_i$, the
two summands of $X$ are $H_1\otimes_{\C} A_2$ and $H_2\otimes_{\C} A_1$.
On the first, $(a_1,a_2)\in D$ acts by $\pi_1(a_1)\otimes 1$ and on
the second by $\pi_2(a_2)\otimes 1$. Thus the left action is faithful
and unital. This reduced free product is therefore locally
reflexive by Proposition~\ref{prop:toeplitz} and
Lemma~\ref{lem:permanence}. We do not use the scalar construction of
\cite[Proposition~4.2]{DS}: it involves $A_1\tmin A_2$, whose local
reflexivity follows from Corollary~\ref{cor:tensor} only when one
factor is exact.

For a finite set $S\subset I$, let $A_S$ denote the reduced free
product of $(A_i,\varphi_i)_{i\in S}$. Each $A_S$ is locally reflexive
by induction on $|S|$. Indeed, for $j\notin S$, associativity
identifies $A_{S\cup\{j\}}$ with the reduced free product of $A_S$
and $A_j$. The free product state on $A_S$ again has faithful GNS
representation, so the two-factor case applies. Now let $A$ be the
reduced free product over $I$. Taking scalar subalgebras outside $S$,
we obtain a canonical embedding of $A_S$ in $A$ by
\cite[Theorem~1.3]{BD}. The subalgebras $A_S$ form a directed family
whose union contains every reduced word. Thus $A$ is their inductive
limit and is locally reflexive by Corollary~\ref{cor:inductive}.

Under the convention allowing arbitrary states, the same argument
applies to the GNS quotients $A_i/\ker\pi_{\varphi_i}$
\cite[Remark~1.4]{BD}, which are locally reflexive by
Lemma~\ref{lem:permanence}.
\end{proof}

\begin{corollary}[Corollary~\ref{cor:kirchberg-fh}(h)]\label{cor:simple}
Every locally reflexive $C^*$-algebra embeds into a simple
locally reflexive $C^*$-algebra.
\end{corollary}

\begin{proof}
Suppose first that $A$ is nonzero and unital. We choose a faithful unital
representation $\pi\colon A\to\BH(K)$. On the right Hilbert
$A$-module $X=K\otimes_{\C}A$, let $A$ act on the left through $\pi$.
This left action is faithful and unital, and $X$ is full since
$\langle\eta\otimes1,\eta\otimes1\rangle=1$ for every unit vector
$\eta\in K$. Thus $X$ satisfies the hypotheses of
Proposition~\ref{prop:toeplitz}.
For $0\ne a\in A_+$, we choose a unit vector $\xi\in K$ with
$\langle\xi,\pi(a)\xi\rangle>0$. As in Kumjian's proof of
\cite[Lemma~2.7]{Kum}, the creation operator
$s_{\xi\otimes1}$ in $\mathcal T(X)$ satisfies
\[
 s_{\xi\otimes1}^*\,a\,s_{\xi\otimes1}
      =\langle\xi,\pi(a)\xi\rangle\,1 .
\]
Thus every nonzero positive element of $A$ generates $\mathcal T(X)$
as a closed ideal.

Let $J$ be a maximal proper closed ideal of the unital algebra
$\mathcal T(X)$. Since $J$ is proper, it contains no nonzero positive
element of $A$. Hence $J\cap A=0$, and the simple quotient
$\mathcal T(X)/J$ contains a faithful copy of $A$. It is locally
reflexive by Proposition~\ref{prop:toeplitz} and
Lemma~\ref{lem:permanence}. For nonunital $A$, apply the unital case to its unitization, which is
locally reflexive by Lemma~\ref{lem:permanence}. The zero algebra
embeds in $\C$.
\end{proof}

\begingroup
\small

\medskip
\noindent
\begin{tabular}{@{}l@{}}
Yafei Zhao\\
Department of Mathematics\\
Zhejiang International Studies University\\
Hangzhou 310023, China\\
E-mail address: yfzhao@zisu.edu.cn\\[1.2ex]
Zhe Dong\\
School of Mathematical Sciences\\
Zhejiang University\\
Hangzhou 310027, China\\
E-mail address: dongzhe@zju.edu.cn\\
Corresponding author.
\end{tabular}
\endgroup

\end{document}